\documentstyle[11pt]{article}
\begin{document}
\author{{Shiqiu Zheng$^1$\thanks{E-mail: shiqiumath@163.com}},\ \ Lidong Zhang$^1$,\ \ Lichao Feng$^2$\\
  \\
\small(1, College of Sciences, Tianjin University of Science and Technology, Tianjin 300457, China)\\ \small(2, College of Sciences, North China University of Science and Technology, Tangshan 063210, China)}
\date{}
\title{\textbf{On the backward stochastic differential equation with generator $f(y)|z|^2$}}\maketitle

\textbf{Abstract:}\quad In this paper, we consider the backward stochastic differential equation (BSDE) with generator $f(y)|z|^2,$ where the function $f$ is defined on an open interval $D$ and locally integrable. The existence and uniqueness of bounded solutions and $L^p(p\geq1)$ solutions of such BSDEs are obtained. Some comparison theorems and a converse comparison theorem of such BSDEs are established. As an application, we give a probabilistic interpretation of viscosity solution of quadratic PDEs.\\



\section{Introduction}
Nonlinear backward stochastic differential equations (BSDEs) were firstly introduced by the seminal paper of
Pardoux and Peng \cite{PP}. Since then, many works have been done to weaken the existence conditions of solutions. An important case is the study of quadratic BSDEs, i.e., the generators have a quadratic growth in control variable $z.$ Quadratic BSDEs have a wide range of applications in stochastic control and finance (see El Karoui and Rouge \cite{El}, Hu, Imkeller and M\"{u}ller \cite{H} and El Karoui, Hamad\`{e}ne and Matoussi \cite{El2}). Using an exponential transform and monotone stability, Kobylanski \cite{K} derived the existence and uniqueness of solutions of quadratic BSDEs with bounded terminal variables. A similar result was obtained by Briand and Elie \cite{BE} using a different method based on Malliavin calculus. Using a localization method, Briand and Hu \cite{BH06}, \cite{BH08} obtained the existence and uniqueness of solutions for exponential integrable terminal variables. Recently, Bahlali et al. \cite{B17} obtained the existence and uniqueness of solutions for $L^2$ integrable terminal variables, but the quadratic term takes the form $f(y)|z|^2,$ where $f$ is globally integrable on $R.$ Yang \cite{Y} further considered the $L^p\ (p\geq1)$ solutions of such quadratic BSDEs, where $f$ is globally integrable on $R$ and bounded on any compact subset of $R.$ A one to one transformation $u_f$ based on $f$ and It\^{o}-Krylov's formula play a crucial role in the works of \cite{B17} and \cite{Y}. Then the following problems arise naturally:
\begin{itemize}
  \item does the BSDE with generator $f(y)|z|^2$ have a solution? when $f$ is defined on an open interval $D$ and locally integrable. If this BSDE has a solution, which space does the solution belong to?
\end{itemize}

The paper is devoted to exploring the above problems. Recently, the usual case that $D=R$ was considered by Bahlali \cite{B19} by a domination method, and a canonical singular case that $f(y)=\frac{1}{y}$ was considered by Bahlali and Tangpi \cite{BT}. In this paper, the function $f$ is defined on an open interval $D$ and locally integrable. As a result, the corresponding BSDEs can be used to deal with some cases that are not covered by \cite{B19} and \cite{BT}. Some instances are given in Example 2.1. Following \cite{B17}, we also use a one to one transformation $u_f$ based on $f$ and It\^{o}-Krylov's formula. To deal with our problems, our transformation $u_f$ is defined on an open interval $D,$ while the $u_f$ in \cite{B19}, \cite{B17} and \cite{Y} are defined on $R$ with $u_f(0)=0.$ The main results of this paper can be considered as a complementary to the corresponding results in \cite{B19}, \cite{B17}, \cite{BT} and \cite{Y}. For example,
\begin{itemize}
 \item The bounded solutions of BSDE$(f(y)|z|^2,\xi)$ are discussed (see Theorem 3.2(i)(ii)). For example, if $f(y)=\frac{-1}{(y-1)(y-6)}, y\in(1,6),$ then for any variable $1<\xi<6,$ such that $u_f(\xi)\in L^1({\cal{F}}_T),$ BSDE$(f(y)|z|^2,\xi)$ will have a solution $(Y_t,Z_t)$ such that $1<Y_t<6.$ Moreover if $2\leq\xi\leq3$, then $2\leq Y_t\leq3.$
 \item Some sufficient conditions are given to guarantee that BSDE$(f(y)|z|^2,\xi)$ will have $L^p(p\geq1)$ solutions (see Theorem 3.2(iii)(iv)(vi)(vii)). For example, if $f(y)=\frac{1}{y^2}+1,$ $y\in(0,\infty),$ then for any positive variable $\xi$ such that $u_f(\xi)\in L^1({\cal{F}}_T),$ BSDE$(f(y)|z|^2,\xi)$ will have a solution $(Y_t,Z_t)\in{\mathcal{S}}^p\times{\cal{H}}^{2p}, (p>1)$. If $f(y)=\frac{|\ln y|}{y},$ $y\in(0,\infty),$ then for any positive variable $\xi$ such that $u_f(\xi)\in L^p({\cal{F}}_T),$ BSDE$(f(y)|z|^2,\xi)$ will have a solution such that $Y_t\in{\mathcal{S}}^p, (p>1)$.
 \item Some (strict) comparison theorems are established for BSDE$(f(y)|z|^2,\xi)$ (see Proposition 4.1 and Proposition 4.3). As an application, we obtain a converse comparison theorem, which shows that the generators can be compared through comparing the solutions of BSDE$(f(y)|z|^2,\xi)$ (see Proposition 4.5).
  \item BSDE$(f(y)|z|^2,\xi)$ is used to give a probabilistic interpretation of viscosity
solution of a quadratic PDE, which contains a coefficient $f$ locally integrable on $D$ (see Theorem 5.1).
 \end{itemize}

This paper is organized as follows. In section 2, we will present some assumptions and lemmas. In Section 3, we will study the existence and uniqueness of bounded solutions and $L^p(p\geq1)$ solutions. In Section 4, we will establish some comparison theorems and a converse comparison theorem. In Section 5, we will give a probabilistic interpretation of viscosity solution of singular quadratic PDEs.
\section{Preliminaries}
Let $(\Omega ,\cal{F},\mathit{P})$ be a complete probability space, ${{(B_t)}_{t\geq
0}}$ be a $d$-dimensional standard Brownian motion on $(\Omega ,\cal{F},\mathit{P})$. Let $({\cal{F}}_t)_{t\geq 0}$ denote
the natural filtration generated by ${{(B_t)}_{t\geq 0}}$, augmented
by the $\mathit{P}$-null sets of ${\cal{F}}$. Let $|z|$ denote its
Euclidean norm, for $\mathit{z}\in R^d$. Let $T>0$ be a given real number and ${\cal{T}}_{0,T}$ be the set of all stopping times $\tau$ satisfying $0\leq \tau\leq T$. Let $D\subset R$ be an open interval which may take one of the following forms $$(a,b), (a,+\infty),(-\infty,b),\ \textmd{or}\ (-\infty,+\infty),$$
where $a$ and $b$ are two real numbers such that $a<b$. We define the following usual spaces

$L_{1,loc}(D)=\{f:\ f:D\rightarrow R,$ is measurable and locally integrable$\}$;

$L_D({\cal{F}}_T)=\{\xi:\ {\cal{F}}_T$-measurable random variable whose range is included in $D\}$;

$L^p({\mathcal {F}}_T)=\{\xi:\ {\cal{F}}_T$-measurable
$R$-valued random variable and ${{E}}\left[|\xi|^p\right]<\infty\},\ p\geq1; $

$L^\infty({\mathcal {F}}_T)=\{\xi:\ {\cal{F}}_T$-measurable
$R$-valued random variable such that $\textrm{esssup}_{\omega\in\Omega}|\xi|<\infty\}; $

${\mathcal{S}}=\{\psi:$ $R$-valued continuous predictable
process$\}$;

${\mathcal{S}}^p=\{\psi:$ process in ${\mathcal{S}}$ such that ${{E}} \left
[{\mathrm{sup}}_{0\leq t\leq T} |\psi _t|^p\right] <\infty \},\ p>0;$

${\mathcal{S}}^\infty=\{\psi:$
process in ${\mathcal{S}}$ and $\|\psi\|_\infty=\textrm{esssup}_{(\omega,t)\in\Omega\times [0,T]}|\psi _t| <\infty\};$

${H}^p=\{\psi:$ $R^d$-valued predictable
process such that $\int_0^T|\psi_t|^pdt
<\infty \},\ p\geq1;$

${\cal{H}}^p=\{\psi:$ $R^d$-valued predictable
process such that ${{E}}\left[\left(\int_0^T|\psi_t|^2dt\right)^{\frac{p}{2}}\right]
<\infty \},\ p>0;$

${\cal{H}}^2_{BMO}=\{\psi:$ process in ${\mathcal{H}}^2$ such that $\sup_{\tau\in{\cal{T}}_{0,T}}E\left[\int_\tau^T|\psi_t|^2dt|{\cal{F}}_\tau\right]<\infty \};$

$W^2_{1,loc}(D)=\{f: \ f:D\rightarrow R,$ such that $f$ and its generalized derivation $f'$ and $f''$ are all locally integrable measurable functions$\}$.\\

For $f\in L_{1,loc}(D),$ we will define a transformation $u^\alpha_f(x)$, which plays an important role in this paper. For $\alpha\in D,$ we define $$u^\alpha_f(x):=\int^x_\alpha\exp\left(2\int^y_\alpha f(z)dz\right)dy,\ \ x\in D.$$\\
\textbf{Remark 2.1}\ \ For $\alpha, \beta\in D$, we have
\begin{eqnarray*}
u^\alpha_f(x)&=&\int^x_\beta\exp\left(2\int^y_\alpha f(z)dz\right)dy+\int^\beta_\alpha\exp\left(2\int^y_\alpha f(z)dz\right)dy\\
&=&\int^x_\beta\exp\left(2\int^y_\beta f(z)dz+2\int^\beta_\alpha f(z)dz\right)dy+\int^\beta_\alpha\exp\left(2\int^y_\alpha f(z)dz\right)dy\\
&=&\exp\left(2\int^\beta_\alpha f(z)dz\right)\int^x_\beta\exp\left(2\int^y_\beta f(z)dz\right)dy+\int^\beta_\alpha\exp\left(2\int^y_\alpha f(z)dz\right)dy,
\end{eqnarray*}
which means that for $\alpha, \beta\in D$, there exist two constants $a>0$ and $b$ such that
$$u^\alpha_f(x)=au^\beta_f(x)+b,\ x\in D.$$

By Remark 2.1, we will find that the different choices of $\alpha\in D$ in the definition of $u^\alpha_f(x)$ do not change the solution of BSDE$(f(y)|z|^2,\xi)$ (see Remark 3.1(i)). Thus, in the following, let $\alpha\in D$ be given and $u^\alpha_f(x)$ be denoted by $u_f(x)$. Let $V:=\{y:y=u_f(x),x\in D\}.$ The following Lemma 2.1 contains some important properties of $u_f(x)$.\\
\\
\textbf{Lemma 2.1}\ \emph{For $f\in L_{1,loc}(D),$ the following properties of $u_f(x)$ hold true.}

\emph{(i) $u_f(x)\in W^2_{1,loc}(D)$, in particular, $u_f(x)\in C^1(D)$;}

\emph{(ii) $u_f(x)$ is strictly increasing;}

\emph{(iii) $u_f''(x)-2f(x)u_f'(x)=0,\ a.e.$ on $D$;}

\emph{(iv) If $g\in L_{1,loc}(D)$ and $g(x)\leq f(x),\ a.e.$ on $D$, then for every $x\in D,$ $u_g(x)\leq u_f(x)$ and for every  $x\in[\alpha,\infty)$, $0<u'_{g}(x)\leq u'_{f}(x)$;}

\emph{(v) $V$ is an open interval and $u_f^{-1}(x)\in W^2_{1,loc}(V)$, in particular, $u_f^{-1}(x)\in C^1(V)$ and is strictly increasing.}\\\\
\emph{Proof.} (i), (iii) and (v) can be derived from Lemma A.1 in \cite{B17}, where the case that $D=R$ is considered. (ii) is clear. Now we prove (iv). Since $g(x)\leq f(x),\ a.e.$ on $D,$ by the monotonicity of integral and exponential function, we can complete this proof.\ \ $\Box$\\

For $f\in L_{1,loc}(D)$ and $\xi\in L_D({\cal{F}}_T),$ this paper considers the following one-dimensional BSDE
$$Y_t=\xi+\int_t^{T}f(Y_s)|Z_s|^2ds-\int_t^{T}Z_sdB_s,\ \ t\in[0,T],\eqno(1)$$
where $\xi$ is the terminal variable and $T$ is the terminal time. (1) is denoted by BSDE$(f(y)|z|^2,\xi)$, has been studied by \cite{B17} and \cite{Y} in the case that $f$ is globally integrable on $R,$ and by \cite{B19} in the case that $f$ is locally integrable on $R.$ In this note, $f$ is assumed to be locally integrable on the open interval $D.$\\\\
\textbf{Definition 2.1} For $f\in L_{1,loc}(D)$ and $\xi\in L_D({\cal{F}}_T),$ the solution of BSDE$(f(y)|z|^2,\xi)$ is a pair $(Y_t,Z_t)\in{\mathcal{S}}
\times{{H}}^2$ satisfying (1) and $\int_0^{T}|f(Y_s)||Z_s|^2ds<\infty,$ and the range of $Y_t$ is included in $D$. \\

The following Example 2.1 shows some common cases contained in our setting, not covered by \cite{B19}, \cite{B17} and \cite{Y}.\\\\
\textbf{Example 2.1}  (i) $f(y)=\frac{|\ln(y)|}{y},$ $y\in D=(0,\infty)$;

(ii)  $f(y)=\frac{-1}{(y-1)(y-6)},$ $y\in D=(1, 6)$;

(iii)  $f(y)=\frac{1}{y^2}+1,$ $y\in D=(0,\infty)$.\\

Using a probabilistic method, the Proposition 2.1 in \cite{B17} established the Krylov's estimate for the solution of BSDE$(f(y)|z|^2,\xi),$ where $f\in L_{1,loc}(R)$. By a slightly modified proof, we obtain the following Krylov's estimate for the case that $f\in L_{1,loc}(D)$. \\\\
\textbf{Lemma 2.2}\ \ \emph{For $f\in L_{1,loc}(D)$ and $\xi\in L_D({\cal{F}}_T),$ let $(Y_t,Z_t)$ be a solution of
BSDE$(f(y)|z|^2,\xi)$, and $K>0$ be a constant such that $[Y_0-K,Y_0+K]\subset D$. Then there exists a positive constant $\gamma$ depending
on $K$ and $\int_{Y_0-K}^{Y_0+K}|f(x)|dx$ such that for any nonnegative measurable
function $\psi\in L_{1,loc}(D)$,}
$$E\int_{0}^{T\wedge\tau_K}\psi(Y_s)|Z_s|^2ds\leq \gamma \int_{Y_0-K}^{Y_0+K}\psi(x)dx,$$
where $\tau_K:=\inf\{t\geq0,Y_t\notin(Y_0-K,Y_0+K)\}.$\\\\
\emph{Proof.}\ \ In fact, if we replace the stopping time $\tau_R$ in the Proposition 2.1 in \cite{B17} by $$\tau_K:=\inf\{t\geq0,Y_t\notin(Y_0-K,Y_0+K)\},$$ then by the same method, we can complete this proof.\ \ $\Box$\\

Using Lemma 2.2 and the same method of Theorem 2.1 in \cite{B17}, we obtain the following It\^{o}-Krylov's formula, which is an extension of It\^{o}'s formula and can be used to deal with BSDE$(f(y)|z|^2,\xi),$ where $f$ is measurable. \\\\
\textbf{Lemma 2.3}\ \ \emph{For $f\in L_{1,loc}(D)$ and $\xi\in L_D({\cal{F}}_T),$ let $(Y_t,Z_t)$ be a solution of
BSDE$(f(y)|z|^2,\xi)$. Then for any $u\in W^2_{1,loc}(D),$ we have }
$$u(Y_t)=u(Y_0)+\int_0^tu'(Y_s)dY_s+\frac{1}{2}\int_0^tu''(Y_s)|Z_s|^2ds, \ \ t\in[0,T].$$
\emph{Proof.}\ \ Let $K>0$ be a constant such that $[Y_0-K,Y_0+K]\subset D$. If we replace the stopping time $\tau_R$ and the interval $[-R,R]$ in the proof of Theorem 2.1 in \cite{B17} by $$\tau_K:=\inf\{t\geq0,Y_t\notin(Y_0-K,Y_0+K)\}$$ and $[Y_0-K,Y_0+K],$ respectively, then by the same method, we can complete this proof.\ \ $\Box$
\section{Existence and uniqueness of BSDE$(f(y)|z|^2,\xi)$}
In this section, we will explore the bounded solutions and $L^p(p\geq1)$ solutions of BSDE$(f(y)|z|^2,\xi),$ where $f$ is defined on an open interval $D$ and locally integrable. For convenience, we note that all the closed subintervals of $D$ or $V$ mentioned in the following will be finite intervals.  We firstly give a necessary condition.\\\\
\textbf{Proposition 3.1}\ \ \emph{Let $f\in L_{1,loc}(D)$ and $\xi\in L_D({\cal{F}}_T).$ If BSDE$(f(y)|z|^2,\xi)$ has a solution $(Y_t,Z_t)$ and there exists a constant $\beta$ such that for every $x\in D$, $u_f(x)\geq\beta$ or for every $x\in D$, $u_f(x)\leq\beta,$ then we have $u_f(\xi)\in L^1({\cal{F}}_T).$}
\\\\
\emph{Proof.} Applying Lemma 2.3 to $u_{f}(Y_t),$ and then by Lemma 2.1(iii), we have
$$u_f(Y_t)=u_f(\xi)+\int_t^Tu'_f(Y_s)Z_sdB_s,\ \ t\in[0,T].\eqno(2)$$
 For $n\geq1,$ we define the following stopping time
$$\tau_n=\inf\left\{t\geq0,\int_0^t|u'_f(Y_s)|^2|Z_s|^2ds\geq n\right\}\wedge T,$$
which implies that $\int_0^{\tau_n\wedge t}u'_f(Y_s)Z_sdB_s$ is a martingale. Then by (2), we have
$$u_f(Y_0)=E[u_f(Y_{\tau_n})].\eqno(3)$$
Since $\tau_n\rightarrow T$, as $n\rightarrow\infty,$ by the continuity of $Y_t$ and $u_f$ (see Lemma 2.1(i)), we have
$$\lim_{n\rightarrow\infty}u_f(Y_{\tau_n})=u_f(\xi).\eqno(4)$$

If there exists a constant $\beta$ such that for every $x\in D$, $u_f(x)\geq\beta,$ then by (4), Fatou's lemma and (3), we deduce
$$0\leq E[u_f(\xi)-\beta]=E[\liminf_{n\rightarrow\infty}(u_f(Y_{\tau_n})-\beta)]\leq \liminf_{n\rightarrow\infty}E[u_f(Y_{\tau_n})-\beta]=u_f(Y_0)-\beta,$$
which implies $u_f(\xi)\in L^1({\cal{F}}_T).$

Similarly, if there exists a constant $\beta$ such that for every $x\in D$, $u_f(x)\leq\beta,$ then by (4), Fatou's lemma and (3), we have
$$0\leq E[\beta-u_f(\xi)]=E[\liminf_{n\rightarrow\infty}(\beta-u_f(Y_{\tau_n}))]\leq \liminf_{n\rightarrow\infty}E[\beta-u_f(Y_{\tau_n})]=\beta-u_f(Y_0),$$
which implies $u_f(\xi)\in L^1({\cal{F}}_T).$\ \ $\Box$\\

By the definition of $u_f$ and Lemma 2.1(iv), it is not hard to check that the following examples satisfy the conditions in Proposition 3.1.\\\\
\textbf{Example 3.1} Let $\beta>0$ be a constant. We have

(i) if $f(y)\geq\beta,\ a.e.$ on $R,$ then for every $x\in R$, $$u_f(x)\geq u_\beta=\frac{1}{2\beta}(\exp(2\beta (x-\alpha))-1)>-\frac{1}{2\beta};$$

(ii) if $f(y)\leq-\beta,\ a.e.$ on $R,$ then for every $x\in R$, $$u_f(x)\leq u_{-\beta}=-\frac{1}{2\beta}(\exp(-2\beta (x-\alpha))-1)\leq\frac{1}{2\beta};$$

(iii) if $f(y)=\frac{\beta}{2y},\ a.e.$ on $(0,\infty),$ then for every $x\in (0,\infty)$, $$u_f(x)=\frac{\alpha}{1+\beta}\left(\left(\frac{x}{\alpha}\right)^{\beta+1}-1\right)\geq-\frac{\alpha}{1+\beta}.$$

The following Theorem 3.2 is the main result of this section.\\\\
\textbf{Theorem 3.2}\ \ \emph{Let $f\in L_{1,loc}(D)$ and $\xi\in L_D({\cal{F}}_T).$ If $u_f(\xi)\in L^1({\cal{F}}_T),$ then BSDE$(f(y)|z|^2,\xi)$ has a unique solution $(Y_t,Z_t)$ such that $u_f(Y_t)=E[u_f(\xi)|{\cal{F}}_t].$ Moreover, we have}

\emph{(i) If $D$ is bounded, then $Y_t\in{\mathcal{S}}^\infty$;}

\emph{(ii) If the range of $\xi$ is included in a closed interval $[a,b]\subset D,$ then $(Y_t,Z_t)\in{\mathcal{S}}^\infty
\times{\cal{H}}^2_{BMO}$ and the range of $Y$ is included in $[a,b]$;}

\emph{(iii) If $\xi\in L^1({\cal{F}}_T),$ $\xi^-\in L^p({\cal{F}}_T)$ and there exists a constant $\beta>0$ such that $f\geq\beta, a.e.$ on $D,$ then $(Y_t,Z_t)\in{\mathcal{S}}^p\times{\cal{H}}^{2p},$ when $p>1,$ and $(Y_t,Z_t)\in{\mathcal{S}}^r\times{\cal{H}}^2, r\in(0,1),$ when $p=1$;}

\emph{(iv) If $\xi\in L^1({\cal{F}}_T),$ $\xi^+\in L^p({\cal{F}}_T)$ and there exists a constant $\beta>0$ such that $f\leq-\beta, a.e.$ on $D,$ then $(Y_t,Z_t)\in{\mathcal{S}}^p\times{\cal{H}}^{2p},$ when $p>1,$ and $(Y_t,Z_t)\in{\mathcal{S}}^r\times{\cal{H}}^2, r\in(0,1),$ when $p=1$;}

\emph{(v) If $u_f(\xi)\in L^p({\cal{F}}_T)$, then $(u_f(Y_t),u'_f(Y_t)Z_t)\in{\mathcal{S}}^p
\times{\cal{H}}^p$, when $p>1,$ and $(u_f(Y_t),u'_f(Y_t)Z_t)\in{\mathcal{S}}^r
\times{\cal{H}}^r, r\in(0,1)$, when $p=1$;}

\emph{(vi) If $\xi\in L^1({\cal{F}}_T),$ $\xi^-\in L^p({\cal{F}}_T)$, $u_f(\xi)\in L^p({\cal{F}}_T)$ and $f\geq0, a.e.$ on $D,$ then $Y_t\in{\mathcal{S}}^p$, when $p>1,$ and $Y_t\in{\mathcal{S}}^r,  r\in(0,1)$, when $p=1$;}

\emph{(vii) If $\xi\in L^1({\cal{F}}_T),$ $\xi^+\in L^p({\cal{F}}_T)$,  $u_f(\xi)\in L^p({\cal{F}}_T)$ and $f\leq0, a.e.$ on $D,$ then $Y_t\in{\mathcal{S}}^p$, when $p>1,$ and $Y_t\in{\mathcal{S}}^r,  r\in(0,1)$, when $p=1$.}\\\\
\emph{Proof.} By Proposition 1.1(i) in \cite{B19}, we get that the BSDE$(0,u_f(\xi))$ has a unique solution $(y_t,z_t)$ such that $y_t=E[u_f(\xi)|{\cal{F}}_t].$ This result can be proved by a general martingale representation theorem (see Corollary 3 in Protter \cite{Pr}, page 189) and a localization method. It follows that the range of $y_t$ is included in $V.$
Thus by Lemma 2.1(v), we can apply Lemma 2.3 to $u_f^{-1}(y_t),$ and then, by setting
$$(Y_t,Z_t):=\left(u_f^{-1}(y_t),\frac{z_t}{u_f'(u_f^{-1}(y_t))}\right),\eqno(5)$$
and Lemma 2.1(iii), we get that the BSDE$(f(y)|z|^2,\xi)$ has a solution $(Y_t,Z_t)$ such that $u_f(Y_t)=E[u_f(\xi)|{\cal{F}}_t].$ Conversely, for a solution $(Y_t,Z_t)$ of BSDE$(f(y)|z|^2,\xi)$ such that $u_f(Y_t)=E[u_f(\xi)|{\cal{F}}_t],$ by Lemma 2.1(i), we can apply Lemma 2.3 to $u_f(Y_t)$, and then by Lemma 2.1(iii), we have
$$u_f(Y_t)=u_f(\xi)+\int_t^Tu'_f(Y_s)Z'_sdB_s,\ \ t\in[0,T],$$ which means that $(u_f(Y_t),u'_f(Y_t)Z_t)$ is a solution of BSDE$(0,u_f(\xi))$ such that $u_f(Y_t)=E[u_f(\xi)|{\cal{F}}_t].$ Then the uniqueness of such solution of BSDE$(f(y)|z|^2,\xi)$ can be obtained from the uniqueness of the solution to BSDE$(0,u_f(\xi))$ and Lemma 2.1(ii).

\textbf{Proof of (i):} Since $D$ is bounded and the range of $Y_t$ is included in $D$, we have $Y_t\in{\mathcal{S}}^\infty$.

\textbf{Proof of (ii):} If the range of $\xi$ is included in a closed interval $[a,b]\subset D,$ then by the fact $y_t=E[u_f(\xi)|{\cal{F}}_t]$ and Lemma 2.1(ii), we get that the range of $y_t$ is included in $[u_f(a),u_f(b)]\subset V.$ This together with Proposition 2.1 in \cite{BE} implies $z_t\in{\cal{H}}^2_{BMO}.$  Then by (5) and Lemma 2.1(v)(i), we get that $(Y_t,Z_t)\in{\mathcal{S}}^\infty
\times{\cal{H}}^2_{BMO}$ and the range of $Y$ is included in $[a,b]$.

\textbf{Proof of (iii): } Since $f>0, a.e.$ on $D,$ by the definition of $u_f$, we get that $u_f'$ is nondecreasing. In view of $(u_f^{-1})'(x)=\frac{1}{u'_f(u_f^{-1}(x))}$ and $u_f^{-1}(x)$ is strictly increasing (see Lemma 2.1(v)),  we get that $(u_f^{-1})'$ is not increasing. This together with the fact $D$ is a convex set implies that $u_f^{-1}$ is concave. Then by (5), the concavity of $u_f^{-1}$ and Jensen's inequality, we have
$$Y_t=u_f^{-1}(y_t)=u_f^{-1}(E[u_f(\xi)|{\cal{F}}_t])\geq E[u_f^{-1}(u_f(\xi))|{\cal{F}}_t]=E[\xi|{\cal{F}}_t]\geq -E[\xi^{-}|{\cal{F}}_t].\eqno(6)$$
Since $\xi^-\in L^p({\cal{F}}_T)$ and $p\geq1,$ $M_t:=-E[\xi^{-}|{\cal{F}}_t]$ is a continuous martingale. By Doob's optional sampling theorem and Jensen's inequality, for $\tau\in{\cal{T}}_{0,T},$ we have for $p\geq1,$
$$E[|M_\tau|^p]=E[|E[\xi^{-}|{\cal{F}}_\tau]|^p]\leq E[E[|\xi^{-}|^p|{\cal{F}}_\tau]]= E[|\xi^-|^p].\eqno(7)$$
For $n\geq1,$ we define the following stopping time
$$\tau_n=\inf\left\{t\geq0,\int_0^tf(Y_s)|Z_s|^2ds\geq n\right\}\wedge T.$$
By the assumption $f\geq\beta, a.e.$ on $D$, (1) and (6), we have
$$ \beta\int_0^{\tau_n}|Z_s|^2ds\leq\int_0^{\tau_n}f(Y_s)|Z_s|^2ds\leq Y_0-M_{\tau_n}+\int_0^{\tau_n}Z_sdB_s.$$
Then by Jensen's inequality, (7) and BDG inequality, we deduce
\begin{eqnarray*}
\ \ \ \ \ \ \beta^p\left(E\left[\left(\int_0^{\tau_n}|Z_s|^2ds\right)^{\frac{p}{2}}\right]\right)^2&\leq& E\left[\left(\beta\int_0^{\tau_n}|Z_s|^2ds\right)^p\right]\\&\leq&3^{p-1}\left(|Y_0|^p+E[|M_{\tau_n}|^p]+E\left(\left|\int_0^{\tau_n}Z_sdB_s\right|^p\right)\right)\\
&\leq&C_p\left(|Y_0|^p+E[|\xi^-|^p]+E\left[\left(\int_0^{\tau_n}|Z_s|^2ds\right)^{\frac{p}{2}}\right]\right), \ \ \ \ \ \ \ \ \ \ (8)
\end{eqnarray*}
where $C_p>0$ is a constant depending only on $p.$ By solving the quadratic inequality (8) with $E[(\int_0^{\tau_n}|Z_s|^2ds)^{\frac{p}{2}}]$ as the unknown variable, we deduce that there exists a constant $K>0$ dependent only on $E[|\xi^-|^p],Y_0,p$ and $\beta,$ such that $E\left[\left(\int_0^{\tau_n}|Z_s|^2ds\right)^{\frac{p}{2}}\right]\leq K.$ Plugging this inequality into (8), we get that there exists a constant $K_1>0$ dependent only on $E[|\xi^-|^p],Y_0,p$ and $\beta,$ such that
$$E\left[\left(\int_0^{\tau_n}|Z_s|^2ds\right)^p\right]\leq K_1.\eqno(9)$$
Since $\tau_n\rightarrow T$, as $n\rightarrow\infty$, by (9) and Fatou's lemma, we have
$$E\left[\left(\int_0^{T}|Z_s|^2ds\right)^p\right]=E\left[\liminf_{n\rightarrow\infty} \left(\int_0^{\tau_n}|Z_s|^2ds\right)^p\right]\leq\liminf_{n\rightarrow\infty} E\left[\left(\int_0^{\tau_n}|Z_s|^2ds\right)^p\right]\leq K_1.\eqno(10)$$
By (1), (6) and the assumption $f\geq\beta, a.e.,D,$ we have
$$M_t\leq Y_{t}\leq Y_0-\int_0^{t}\beta|Z_s|^2ds+\int_0^{t}Z_sdB_s.\eqno(11)$$
Clearly, $|M_t|$ is a submartingale. Then by Doob's maximal inequality (see Karatzas and Shreve \cite{ks}, page 14), we have for $p>1,$
$$E\left[\sup_{t\in[0,T]}[|M_t|^p]\right]\leq \left(\frac{p}{p-1}\right)E[|\xi^-|^p].\eqno(12)$$
Then by (11), (12), BDG inequality and (10), we deduce that, when $p>1$,
\begin{eqnarray*}
\ \ &&E\left[\sup_{t\in[0,T]}|Y_{t}|^p\right]\\&\leq& E\sup_{t\in[0,T]}\left[4^{p-1}\left (|M_t|^p+|Y_0|^p+\left|\int_0^{t}\beta|Z_s|^2ds\right|^p+\left|\int_0^{t}Z_sdB_s\right|^p\right)\right] \\&\leq& 4^{p-1}\left(E\left[\sup_{t\in[0,T]}[|M_t|^p]\right]+|Y_0|^p+E\left[\left(\beta\int_0^T|Z_s|^2ds\right)^p\right]+E\left[\sup_{t\in[0,T]}\left|\int_0^{t}Z_sdB_s\right|^p\right]\right)\\&\leq& K_2\left(E[|\xi^-|^p]+|Y_0|^p+E\left[\left(\int_0^T|Z_s|^2ds\right)^p\right]+E\left[\left(\int_0^{T}|Z_s|^2ds\right)^{\frac{p}{2}}\right]\right)\\
&<&\infty,\ \ \ \ \  \ \ \ \  \ \ \ \  \ \ \ \  \ \ \ \  \ \ \ \  \ \ \ \  \ \ \ \  \ \ \ \  \ \ \ \  \ \ \ \  \ \ \ \  \ \ \ \  \ \ \ \  \ \ \ \  \ \ \ \  \ \ \ \  \ \ \ \  \ \ \ \  \ \ \ \  \ \ \ \  \ \ \ \  \ \ \ \  \ \ \ \  \ \ \ \  \ \ \ \ (13)
\end{eqnarray*}
where $K_2>0$ is a constant dependent only on $p$ and $\beta$. Then by (10) and (13), we have $(Y_t,Z_t)\in{\mathcal{S}}^p\times{\cal{H}}^{2p},$ when $p>1.$

Now, we consider the case: $p=1.$ By Lemma 6.1 in Briand et al. \cite{lp}, we have for $r\in(0,1),$
$$E\left[\sup_{t\in[0,T]}|M_t|^r\right]\leq \frac{1}{1-r}[E|\xi^-|]^r.\eqno(14)$$
When $p=1,$ by (11), (14), BDG inequality and (10), we deduce for $r\in(0,1),$
\begin{eqnarray*}
\ \ \ \ &&E\left[\sup_{t\in[0,T]}|Y_{t}|^r\right]\\&\leq& E\sup_{t\in[0,T]}\left[\left (|M_t|^r+|Y_0|^r+\left|\int_0^{t}\beta|Z_s|^2ds\right|^r+\left|\int_0^{t}Z_sdB_s\right|^r\right)\right] \\&\leq&
E\left[\sup_{t\in[0,T]}|M_t|^r\right]+|Y_0|^r+E\left[\left(\int_0^{T}\beta|Z_s|^2ds\right)^r\right]+E\left(\sup_{t\in[0,T]}\left|\int_0^{t}Z_sdB_s\right|^r\right)\\&\leq& K_3\left([E|\xi^-|]^r+|Y_0|^r+E\left[\left(\int_0^{T}|Z_s|^2ds\right)^r\right]+E\left[\left(\int_0^{T}|Z_s|^2ds\right)^{\frac{r}{2}}\right]\right)\\
&<&\infty,\ \ \ \ \  \ \ \ \  \ \ \ \  \ \ \ \  \ \ \ \  \ \ \ \  \ \ \ \  \ \ \ \  \ \ \ \  \ \ \ \  \ \ \ \  \ \ \ \  \ \ \ \  \ \ \ \  \ \ \ \  \ \ \ \  \ \ \ \ \ \  \ \ \ \ \ \ \ \ \  \ \ \ \  \ \ \ \  \ \ \ \  \ \ \ \  \ \ \ \ \ \ \ (15)
\end{eqnarray*}
where $K_3>0$ is a constant depending only on $r$ and $\beta$. Then by (10) and (15), we have $(Y_t,Z_t)\in{\mathcal{S}}^r\times{\cal{H}}^2, r\in(0,1),$ when $p=1.$

\textbf{Proof of (iv):}   The proof is similar as (iii). Since $f<0, a.e.$ on $D,$ by the definition of $u_f$, we get that $u_f'$ is not increasing. In view of $(u_f^{-1})'(x)=\frac{1}{u'_f(u_f^{-1}(x))}$ and $u_f^{-1}(x)$ is strictly increasing,  we get that $(u_f^{-1})'$ is not decreasing. This together with the fact $D$ is a convex set implies that $u_f^{-1}$ is convex. Then by (5), the convexity of $u_f^{-1}$ and Jensen's inequality, we have
$$Y_t=u_f^{-1}(y_t)=u_f^{-1}(E[u_f(\xi)|{\cal{F}}_t])\leq E[u_f^{-1}(u_f(\xi))|{\cal{F}}_t]= E[\xi|{\cal{F}}_t]\leq E[\xi^{+}|{\cal{F}}_t].\eqno(16)$$
Since $\xi^+\in L^p({\cal{F}}_T)$ and $p\geq1$, we get that $N_t:=E[\xi^{+}|{\cal{F}}_t]$ is a continuous martingale. Then by Doob's optional sampling theorem and Jensen's inequality, for $\tau\in{\cal{T}}_{0,T},$ we have for $p\geq1$
$$E[|N_\tau|^p]=E[|E[\xi^{+}|{\cal{F}}_\tau]|^p]\leq E[E[|\xi^{+}|^p|{\cal{F}}_\tau]]=E[|\xi^+|^p].\eqno(17)$$
For $n\geq1,$ we define the following stopping time
$$\tau_n=\inf\left\{t\geq0,\int_0^t-f(Y_s)|Z_s|^2ds\geq n\right\}\wedge T.$$
By the assumption $f\leq-\beta, a.e.$ on $D$, (1) and (16), we have
$$ -\beta\int_0^{\tau_n}|Z_s|^2ds\geq\int_0^{\tau_n}f(Y_s)|Z_s|^2ds\geq Y_0-N_{\tau_n}+\int_0^{\tau_n}Z_sdB_s.$$
Then by Jensen's inequality, (17) and BDG inequality, we deduce
\begin{eqnarray*}
\ \ \ \ \ \ \ \beta^p\left(E\left[\left(\int_0^{\tau_n}|Z_s|^2ds\right)^{\frac{p}{2}}\right]\right)^2&\leq& E\left[\left(\beta\int_0^{\tau_n}|Z_s|^2ds\right)^p\right]\\&\leq& 3^{p-1}\left(|Y_0|^p+E[|N_{\tau_n}|^p]+E\left(\left|\int_0^{\tau_n}Z_sdB_s\right|^p\right)\right)\\
&\leq&c_p\left(|Y_0|^p+E[|\xi^+|^p]+E\left[\left(\int_0^{\tau_n}|Z_s|^2ds\right)^{\frac{p}{2}}\right]\right), \ \ \ \ \ \ \ (18)
\end{eqnarray*}
where $c_p>0$ is a constant depending only on $p.$ Then by solving the quadratic inequality (18) with $E[(\int_0^{\tau_n}|Z_s|^2ds)^{\frac{p}{2}}]$ as the unknown variable, we deduce that there exists a constant $K_4>0$ dependent only on $E[|\xi^+|^p],Y_0,p$ and $\beta,$ such that $E[(\int_0^{\tau_n}|Z_s|^2ds)^{\frac{p}{2}}]\leq K_4.$ Then by a similar argument as (10), we deduce
$$E\left[\left(\int_0^{T}|Z_s|^2ds\right)^p\right]<\infty.\eqno(19)$$
By the assumption $f\leq -\beta, a.e.$ on $D,$ (1) and (16), we have
$$N_t\geq Y_{t}\geq Y_0+\int_0^{t}\beta|Z_s|^2ds+\int_0^{t}Z_sdB_s.\eqno(20)$$
In view of (16)-(20), then by similar arguments as (13) and (15), we can complete this proof.

\textbf{Proof of (v):} By Theorem 4.2 and Theorem 6.3 in \cite{lp}, we have $(y_t,z_t)\in{\mathcal{S}}^p\times{\cal{H}}^p,$ when $p>1$ and $(y_t,z_t)\in{\mathcal{S}}^r\times{\cal{H}}^r, r\in(0,1),$ when $p=1.$ Then by (5), we can complete this proof.

\textbf{Proof of (vi):} Since $f\geq0, a.e.$ on $D,$ then by a similar argument as (6), we have $$Y_t\geq -E[\xi^{-}|{\cal{F}}_t].\eqno(21)$$
By Lemma 2.1(iv), we have $u_f(x)\geq u_0(x)=x-\alpha.$ In view of $u_f(x)$ and $u_0(x)$ are both increasing in $x$, we have $u_f^{-1}(x)\leq u_0^{-1}(x)=x+\alpha.$ This together with (5) implies
$$Y_t=u_f^{-1}(y_t)\leq y_t+\alpha.\eqno(22)$$
Since $\xi^-\in L^p({\cal{F}}_T),$ by (12) and (14), we get $-E[\xi^{-}|{\cal{F}}_t]\in{\mathcal{S}}^p,$ when $p>1$ and $-E[\xi^{-}|{\cal{F}}_t]\in{\mathcal{S}}^r, r\in(0,1),$ when $p=1.$ On the other hand, from the proof of (v), we get $(y_t,z_t)\in{\mathcal{S}}^p\times{\cal{H}}^p,$ when $p>1$ and $(y_t,z_t)\in{\mathcal{S}}^r\times{\cal{H}}^r, r\in(0,1),$ when $p=1$. Then by (21) and (22), we obtain (vi).

\textbf{Proof of (vii):}  The proof is similar as (vi). Since $f\leq0, a.e.$ on $D,$ then by a similar argument as (16), we have $$Y_t\leq E[\xi^{+}|{\cal{F}}_t].\eqno(23)$$
By Lemma 2.1(iv), we have $u_f(x)\leq u_0(x)=x-\alpha.$ In view of $u_f(x)$ and $u_0(x)$ are both increasing in $x$, we have $u_f^{-1}(x)\geq u_0^{-1}(x)=x+\alpha.$ This together with (5) implies
$$Y_t=u_f^{-1}(y_t)\geq y_t+\alpha.\eqno(24)$$
Then by (23), (24), and similar arguments as (vi), we obtain (vii). \ \ $\Box$\\\\
\textbf{Remark 3.1}
\begin{itemize}
\item (i)\ \ In Theorem 3.2, in view of the fact that $u_f(Y_t)=E[u_f(\xi)|{\cal{F}}_t]$ and Remark 2.1, we can check that the solution $(Y_t,Z_t)$ will not change for a different choice of $\alpha\in D$ in the definition of $u_f(x)$.
\item (ii)\ \ In Theorem 3.2, if $u_f(\xi)\in L^\infty({\cal{F}}_T)$, does $Y_t\in{\mathcal{S}}^\infty$? In fact, the result is not always true. For example, let $D=R,\ f=\frac{1}{2}$ and $\xi$ be an unbounded negative random variable. We choose $\alpha=0.$ Clearly, $u_f(\xi)=\exp(\xi)-1$ is bounded. But if $Y_t\in{\mathcal{S}}^\infty$, then $\xi$ will be bounded. This induces a contradiction.
\item (iii)\ \ In Theorem 3.2(ii), if the range of $\xi$ is included in $[a,b]\cup[c,d]\subset D,$ where $d>c>b>a$, does the range of $Y$ is included in $[a,b]\cup[c,d]$? In fact, the result is not always true. For example, let $f(y)=\frac{1}{2y},\ y\in(0,10), a=1, b=3, c=4, d=6,$ and  $P(\xi=2)=P(\xi=5)=\frac{1}{2}.$ We choose $\alpha=1,$ then $u_f(x)=\frac{x^2}{2}-\frac{1}{2}.$ Clearly, the BSDE$(0,u_f(\xi))$ has a unique solution $(y_t,z_t)\in{\mathcal{S}}^\infty\times{\cal{H}}^2_{BMO}$ such that $y_0=E[u_f(\xi)]=\frac{27}{4}.$ But if the range of $Y$ is included in $[1,3]\cup[4,6],$ then by (5), we will have $y_0=u_f(Y_0)\in[0,4]\cup[\frac{15}{2},\frac{35}{2}].$ This induces a contradiction.
\item (iv)\ \ In Theorem 3.2(iii), if $f(x)=0,x\in D,$ then $u_f(x)=x-\alpha,$ and the corresponding BSDE$(f(y)|z|^2,\xi)$ becomes
    $$Y_t=\xi-\int_t^TZ_sdB_s$$ with $\xi\in L^1({\cal{F}}_T)$ and $\xi^-\in L^p({\cal{F}}_T)$. Clearly, if $\xi\notin L^p({\cal{F}}_T)$, then we will not have $Y_t\in{\mathcal{S}}^p.$ Similarly, in Theorem 3.2(iv), if $f(x)=0,x\in D$, then we will also not always have $Y_t\in{\mathcal{S}}^p.$ As a result, we must assume $\beta>0$ in Theorem 3.2(iii)(iv). This strengthens the integrability of terminal variables. By Example 3.1(i), in Theorem 3.2(iii), we have $u_f(x)\geq\frac{1}{2\beta}(\exp(2\beta (x-\alpha))-1),$ which implies that the integrability of $\xi$ is not weaker than $E[\exp(2\beta\xi)]<\infty.$ By Example 3.1(ii), in Theorem 3.2(iv), we have $u_f(x)\leq-\frac{1}{2\beta}(\exp(-2\beta (x-\alpha))-1),$ which implies that the integrability of $\xi$ is not weaker than $E[\exp(-2\beta\xi)]<\infty.$
\item (v)\ \ Comparing with 3.2(iii)-(vii), Theorem 4 in \cite{Y} considered $L^p$-solutions when $\xi\in L^p({\cal{F}}_T)$ and $f$ is integrable on $R$. Proposition 2.2 in \cite{B19} gave a condition which guarantees $Z_t\in{\cal{H}}^2,$ when $u_f(\xi)\in L^1({\cal{F}}_T)$ and $f$ is locally integrable on $R$. Proposition 3.2 in \cite{B19} gave a condition which guarantees $u_f(Y_t)\in{\mathcal{S}}^p$ and $Z_t\in{\cal{H}}^2,$ when $u_f(\xi)\in L^p({\cal{F}}_T)$ and $f$ is locally integrable on $R$. Proposition 3.4 in \cite{BT} gave a condition which guarantees $(Y_t,Z_t)\in{\mathcal{S}}^2\times{\cal{H}}^2,$ when $\xi\in L^{2\delta+1}({\cal{F}}_T)$ and $f(y)=\frac{\delta}{y}, \delta\neq\frac{1}{2}$.
\end{itemize}

In the following Proposition 3.3, we consider a slight generalization of BSDE$(f(y)|z|^2,\xi)$. \\\\
\textbf{Proposition 3.3}\ \ \emph{Let $f\in L_{1,loc}(D)$ and $K\in R.$ If $\xi\in L_D({\cal{F}}_T)$ and $u_f(\xi)\in L^p({\cal{F}}_T), p>1$, then the BSDE$(f(y)|z|^2+K|z|,\xi)$ has a unique solution $(Y_t,Z_t)$ such that $(u_f(Y_t),u'_f(Y_t)Z_t)\in{\mathcal{S}}^p
\times{\cal{H}}^p$. Moreover, if the range of $\xi$ is included in a closed interval $[a,b]\subset D,$ then $(Y_t,Z_t)$ is unique in ${\mathcal{S}}^\infty\times{\cal{H}}^2_{BMO},$ and the range of $Y$ is included in $[a,b]$.}
\\\\
\emph{Proof.}\ \ Since $u_f(\xi)\in L^p({\cal{F}}_T), p>1$, by Theorem 4.2 in \cite{lp}, we get that the following BSDE $$y_t=u_f(\xi)+\int_t^{T}K|z_s|ds-\int_t^{T}z_sdB_s,\ \ t\in[0,T],\eqno(25)$$ has a unique solution $(y_t,z_t)\in{\mathcal{S}}^p
\times{\cal{H}}^p$. Set $\theta_s:=\frac{K}{|z_s|}z_s1_{\{|z_s|>0\}}, s\in[0,T].$ Then BSDE(25) can be rewritten as
$$y_t=\xi+\int_t^{T}\theta_sz_sds-\int_t^{T}z_sdB_s,\ \ t\in[0,T].$$
Since $|\theta_t|\leq K,$ by Theorem 4.2 in \cite{lp} again, we get that $(y_t,z_t)$ is the unique $L^p$ solution of the above linear BSDE. By Theorem 4.4 on page 61 and Lemma 2.3 on page 93 in Mao \cite{M}, we get that the following linear SDE
$$X_t=1+\int_0^t\theta_sX_sdB_s,\ \ t\in[0,T],$$
has a unique solution $X_t\in{\mathcal{S}}^{\frac{p}{p-1}}$ with
$$X_t=\exp\left\{\int_0^t\theta_sdB_s-\frac{1}{2}\int_0^t|\theta_s|^2ds\right\},\ \ s\in[t,T].\eqno(26)$$ Applying It\^{o}'s formula to $X_ty_t$, we can show that $X_ty_t$ is a local martingale. Since $y_t\in{\mathcal{S}}^p$ and $X_t\in{\mathcal{S}}^{\frac{p}{p-1}}$, by H\"{o}lder's inequality, we have $E[\sup_{t\in[0,T]}|X_ty_t|]<\infty,$ which implies that $X_ty_t$ is a uniformly integrable martingale.
Then for $t\in[0,T],$ we have $$X_ty_t=E[X_Tu_f(\xi)|{\cal{F}}_t].$$
Then by (26) and Girsanov's theorem, we get $$y_t=E\left[\frac{X_T}{X_t}u_f(\xi)|{\cal{F}}_t\right]=E_Q[u_f(\xi)|{\cal{F}}_t],\eqno(27)$$ where $Q$ is the probability measure satisfying $\frac{dQ}{dP}=X_T.$
By (27), we get that the range of $y_t$ is included in $V.$
Thus by Lemma 2.1(v), we can apply Lemma 2.3 to $u_f^{-1}(y_t),$ then by setting
$$(Y_t,Z_t):=\left(u_f^{-1}(y_t),\frac{z_t}{u_f'(u_f^{-1}(y_t))}\right),\eqno(28)$$
and Lemma 2.1(iii), we have
$$Y_t=\xi+\int_t^{T}(f(Y_s)|Z_s|^2+K|Z_s|)ds-\int_t^{T}Z_sdB_s,\ \ t\in[0,T],$$
which means that BSDE$(f(y)|z|^2+K|z|,\xi)$ has a solution $(Y_t,Z_t)$ such that $(u_f(Y_t),u'_f(Y_t)Z_t)\in{\mathcal{S}}^p
\times{\cal{H}}^p$. Conversely, for a solution $(Y_t,Z_t)$ of BSDE$(f(y)|z|^2+K|z|,\xi)$ such that $(u_f(Y_t),u'_f(Y_t)Z_t)\in{\mathcal{S}}^p
\times{\cal{H}}^p$, by Lemma 2.1(i), we can apply Lemma 2.3 to $u_f(Y_t)$, and then by Lemma 2.1(iii), we get that $(u_f(Y_t),u'_f(Y_t)Z_t)$ is a solution of BSDE$(K|z|,u_f(\xi))$. Then the uniqueness of solution of BSDE$(f(y)|z|^2+K|z|,\xi)$ can be obtained from the uniqueness of solution of BSDE$(K|z|,u_f(\xi))$ and Lemma 2.1(ii).

Moreover, if the range of $\xi$ is included in a closed interval $[a,b]\subset D,$ then by Lemma 2.1(ii), we get that the range of $u_f(\xi)$ is included in $[u_f(a),u_f(b)].$ Then by Theorem 2.1 and Proposition 2.1 in \cite{BE}, we have $(y_t,z_t)\in{\mathcal{S}}^\infty
\times{\cal{H}}^2_{BMO}$. Moreover, by (27), we get that the range of $y_t$ is included in $[u_f(a),u_f(b)].$ Then by (28) and Lemma 2.1(v), we get that the range of $Y$ is included in $[a,b].$\ \ $\Box$
\section{Comparison and converse comparison of BSDE$(f(y)|z|^2,\xi)$}
In this section, we will study the comparison and converse comparison theorems for BSDEs of type BSDE$(f(y)|z|^2,\xi),$ and give some applications. The following Proposition 4.1 is a (strict) comparison theorem, whose proof follows the method of Proposition 3.2 in \cite{B17}.    \\
\\
\textbf{Proposition 4.1}\ \ \emph{Let $f_1\in L_{1,loc}(D), K\geq0$, $g(t)\in H^1,$ $\xi_i\in L_D({\cal{F}}_T)$ and $u_{f_1}(\xi_i)\in L^1({\cal{F}}_T),$ $ i=1,2.$ If $\xi_1\leq\xi_2,$ and BSDE$(f_1(y)|z|^2-K|z|,\xi_1)$ and BSDE$(g(t),\xi_2)$ have the solutions $(Y^1_t,Z^1_t)$ and $(Y^2_t,Z^2_t),$ respectively, such that $(\int_0^tu'_{f_1}(Y^i_s)Z^i_sdB_s)_{t\in[0,T]}$ is a martingale, and for each $t\in[0,T]$, $u_{f_1}(Y^i_t)\in L^1({\cal{F}}_T)$ and $f_1(Y_t^2)|Z_t^2|^2\leq g(t)$, $i=1,2,$ then for each $t\in[0,T],$ we have $Y^1_t\leq Y^2_t.$ Moreover, if $P(\xi_1<\xi_2)>0$, then for each $t\in[0,T],$ we have $P(Y^1_t<Y^2_t)>0,$ or if $\lambda\otimes P\left(f_1(Y_t^2)|Z_t^2|^2<g(t)\right)>0$\footnote{$\lambda$ is the Lebesgue measure and $\lambda\otimes P$ is the product measure of $\lambda$ and $P$}, then we have $Y^1_0<Y^2_0.$} \\\\
\emph{Proof.}  Applying Lemma 2.3 to $u_{f_1}(Y^1_t),$ and then by Lemma 2.1(iii), we have
$$u_{f_1}(Y^1_t)=u_{f_1}(\xi_1)-\int_t^TK|Z^1_s|ds-\int_t^Tu'_{f_1}(Y^1_s)Z^1_sdB_s.\eqno(29)$$
Applying Lemma 2.3 to $u_{f_1}(Y^2_t),$ and by the assumption $f_1(Y_t^2)|Z_t^2|\leq g(t)$ and Lemma 2.1(iii), we have
\begin{eqnarray*}
\ \ u_{f_1}(Y^2_t)&=&u_{f_1}(\xi_2)+\int_t^Tu'_{f_1}(Y^2_s)g(s)ds-\int_t^T\frac{1}{2}u''_{f_1}(Y^2_s)|Z^2_s|^2ds-\int_t^Tu'_{f_1}(Y^2_s)Z^2_sdB_s\\&\geq&u_{f_1}(\xi_2)+\int_t^Tf_1(Y^2_s)u'_{f_1}(Y^2_s)|Z^2_s|^2ds
 -\int_t^T\frac{1}{2}u''_{f_1}(Y^2_s)|Z^2_s|^2ds-\int_t^Tu'_{f_1}(Y^2_s)Z^2_sdB_s\\
&=& u_{f_1}(\xi_2)-\int_t^Tu'_{f_1}(Y^2_s)Z^2_sdB_s.\ \ \ \ \ \ \ \ \ \ \  \ \ \ \ \  \ \ \ \ \ \ \ \ \ \ \ \ \ \  \ \ \ \ \ \ \ \ \ \ \  \ \ \ \ \ \ \ \ \ \  \ \ \ \ \ \ \ \ \ \  \ (30)
\end{eqnarray*}
By the assumptions of the proposition, (29), (30) and Lemma 2.1(ii),  we get for each $t\in[0,T],$
$$u_{f_1}(Y^2_t)\geq E[u_{f_1}(\xi_2)|{\cal{F}}_t]\geq E[u_{f_1}(\xi_1)|{\cal{F}}_t]\geq u_{f_1}(Y^1_t).
\eqno(31)$$
Taking the inverse transformation $u^{-1}_{f_1}$ to (31), we have $Y^1_t\leq Y^2_t.$
Moveover, if $P(\xi_1<\xi_2)>0$, then by (31) and the fact that $u_{f_1}$ and $u^{-1}_{f_1}$ are both strictly increasing, we have for each $t\in[0,T],$ $P(Y^1_t<Y^2_t)>0.$ If $\lambda\otimes P\left(f_1(Y_t^2)|Z_t^2|<g(t)\right)>0$, then by (30) and Lemma 2.1(ii), we get that $u_{f_1}(Y^2_0)>E[u_{f_1}(\xi_2)]$. By this inequality and (31), we have $u_{f_1}(Y^2_0)>u_{f_1}(Y^1_0)$.
Hence, by the fact that $u^{-1}_{f_1}$ is strictly increasing, we have $Y^1_0<Y^2_0.$\ $\Box$\\

By a similar argument, we have\\\\
\textbf{Corollary 4.2}\ \emph{Let $f_2\in L_{1,loc}(D), K\geq0$, $g(t)\in H^1,$ $\xi_i\in L_D({\cal{F}}_T)$ and $u_{f_2}(\xi_i)\in L^1({\cal{F}}_T),$ $i=1,2$. If  $\xi_1\leq\xi_2$, and BSDE$(g(t),\xi_1)$ and BSDE$(f_2(y)|z|^2+K|z|,\xi_2)$ have the solutions $(Y^1_t,Z^1_t)$ and $(Y^2_t,Z^2_t),$ respectively, such that $(\int_0^tu'_{f_2}(Y^i_s)Z^i_sdB_s)_{t\in[0,T]}$ is a martingale, and for each $t\in[0,T]$, $u_{f_2}(Y^i_t)\in L^1({\cal{F}}_T)$ and $g(t)\leq f_2(Y_t^1)|Z_t^1|^2, i=1,2,$ then for each $t\in[0,T],$ we have $Y^1_t\leq Y^2_t.$ Moreover, if $P(\xi_1<\xi_2)>0$, then for each $t\in[0,T],$ we have $P(Y^1_t<Y^2_t)>0,$ or if $\lambda\otimes P\left(g(t)<f_2(Y_t^1)|Z_t^1|^2\right)>0$, then we have $Y^1_0<Y^2_0.$} \\\\
\textbf{Proposition 4.3}\ \ \emph{Let $f_i\in L_{1,loc}(D)$, $\xi_i\in L_D({\cal{F}}_T)$ and $u_{f_i}(\xi_i)\in L^p({\cal{F}}_T), p>1, i=1,2.$ Let $(Y^i_t,Z^i_t)$ be the solution of BSDE$(f_i(y)|z|^2,\xi_i)$, such that $u_{f_i}(Y_t^i)=E[u_{f_i}(\xi_i)|{\cal{F}}_t], i=1,2$. If $\xi_1\leq\xi_2$, $f_1\leq f_2$, and one of the following two conditions holds}

\emph{(a1) there exists a constant $\zeta\in D$ such that $\xi_2\geq\zeta,$ }

\emph{(a2) there exists a constant $\beta$ such that $u_{f_1}(x)\geq\beta,$}\\
\emph{then for each $t\in[0,T],$ we have $Y^1_t\leq Y^2_t.$ Moreover, if $P(\xi_1<\xi_2)>0$, then for each $t\in[0,T],$ we have $P(Y^1_t<Y^2_t)>0.$ In particular, $Y^1_0<Y^2_0.$}  \\\\
\emph{Proof.} We only consider the case that $D=(a,\infty).$ The arguments of the other cases are similar.
For $n\geq1,$ let $a_n=a+\frac{1}{n}$ and $u_{f_{i,n}}(\cdot)$ be the transformation $u_{f_{i}}(\cdot)$ defined with $u_{f_{i}}(a_n)=0,$ $i=1,2.$ By Theorem 3.2(v) and Remark 3.1(i), we have $(u_{f_{i,n}}(Y^i_t),u'_{f_{i,n}}(Y^i_t)Z^i_t)\in{\mathcal{S}}^p
\times{\cal{H}}^p.$

(a1)  By Proposition 4.1, it is enough to prove that $(\int_0^tu'_{f_{1,n}}(Y^2_s)Z^2_sdB_s)_{t\in[0,T]}$ is a martingale and $u_{f_{1,n}}(Y^2_t)\in{\mathcal{S}}^p.$ In fact, if there exists a constant $\zeta\in D$ such that $\xi_2\geq\zeta,$ then by Lemma 2.1(ii)(v), we have
$$Y^2_t=u_{f_2}^{-1}(E[u_{f_2}(\xi_2)|{\cal{F}}_t])\geq u_{f_2}^{-1}(E[u_{f_2}(\zeta)|{\cal{F}}_t])=\zeta.$$
Thus, there exists a constant $K>0$ such that for any $n>K$, $Y^2_t\geq a_n.$ Let $n>K.$ By Lemma 2.1(iv), we have
$$\int_0^T|u'_{f_{1,n}}(Y^2_s)|^2|Z^2_s|^2ds\leq \int_0^T|u'_{f_{2,n}}(Y^2_s)|^2|Z^2_s|^2ds,$$
which, together with the fact $u'_{f_{2,n}}(Y^2_t)Z^2_t\in{\cal{H}}^p$, implies that $(\int_0^tu'_{f_{1,n}}(Y^2_s)Z^2_sdB_s)_{t\in[0,T]}$ is a uniformly integrable martingale. By Lemma 2.1(iv), we have $u_{f_{1,n}}(\zeta)\leq u_{f_{1,n}}(Y^2_t)\leq u_{f_{2,n}}(Y^2_t),$ which together with the fact $u_{f_{2,n}}(Y^2_t)\in{\mathcal{S}}^p
$ implies $u_{f_{1,n}}(Y^2_t)\in{\mathcal{S}}^p.$ The proof is complete.

(a2)  We define the following stopping time
$$\tau_n=\inf\{t\geq0, Y^2_t\leq a_n\}\wedge T.$$
Then there exists a constant $K>0$ such that for any $n>K$, we have $\tau_n>0.$ Let $n>K.$ By Lemma 2.1(iv), we have
$$\int_0^{\tau_n}|u'_{f_{1,n}}(Y^2_s)|^2|Z^2_s|^2ds\leq \int_0^{\tau_n}|u'_{f_{2,n}}(Y^2_s)|^2|Z^2_s|^2ds.$$
which, together with the fact $u'_{f_{2,n}}(Y^2_t)Z^2_t\in{\cal{H}}^p$, implies that $(\int_0^{\tau_n\wedge t}u'_{f_{1,n}}(Y^2_s)Z^2_sdB_s)_{t\in[0,T]}$ is a martingale. Since $u_{f_1}(x)\geq\beta, x\in D,$ we have $\beta\leq u_{f_{1,n}}(Y^2_t)\leq u_{f_{2,n}}(Y^2_t),$ which together with the fact $u_{f_{2,n}}(Y^2_t)\in{\mathcal{S}}^p
$ implies $u_{f_{1,n}}(Y^2_t)\in{\mathcal{S}}^p.$  Then by a similar argument as (30) (replace $u_{f_{1}}$ and $g(t)$ in (30) by $u_{f_{1,n}}(\cdot)$ and $f_2(Y_t^2)|Z_t^2|^2$, respectively), we deduce
$$u_{f_{1,n}}(Y^2_{\tau_n\wedge t})\geq E[u_{f_{1,n}}(Y^2_{\tau_n})|{\cal{F}}_t].$$
Thanks to Remark 2.1, we have
$$u_{f_{1}}(Y^2_{\tau_n\wedge t})\geq E[u_{f_{1}}(Y^2_{\tau_n})|{\cal{F}}_t].$$
Since $\tau_n\rightarrow T,$ as $n\rightarrow\infty,$ passing to the limit in the inequality above, we have
$$u_{f_{1}}(Y^2_{t})\geq \liminf_{n\rightarrow\infty}E[u_{f_{1}}(Y^2_{\tau_n})|{\cal{F}}_t].$$
Since $u_{f_1}(x)\geq\beta, x\in D,$ then by Fatou's Lemma, we deduce
$$u_{f_{1}}(Y^2_{t})-\beta\geq \liminf_{n\rightarrow\infty}E[u_{f_{1}}(Y^2_{\tau_n})-\beta|{\cal{F}}_t]\geq E[u_{f_{1}}(\xi_2)-\beta|{\cal{F}}_t]\geq E[u_{f_1}(\xi_1)-\beta|{\cal{F}}_t]=u_{f_1}(Y^1_t)-\beta.$$
which means that
$$u_{f_{1}}(Y^2_{t})\geq E[u_{f_{1}}(\xi_2)|{\cal{F}}_t]\geq E[u_{f_1}(\xi_1)|{\cal{F}}_t]=u_{f_1}(Y^1_t).$$
Then by similar arguments as in the proof of Proposition 4.1 above, we can complete this proof. \ \ $\Box$\\

Now, we will present two applications of comparison theorems. One is Proposition 4.4 which is an existence result of bounded solutions, the other is Proposition 4.5 which is a converse comparison theorem. \\\\
\textbf{Assumption (A)}
(i) $F_D(s,\omega,y,z): [0,T]\times \Omega\times D\times R^d\mapsto R,$ is continuous in $(y,z)$ for any $(s,\omega)\in[0,T]\times \Omega$;

(ii)  There exist a positive, locally bounded $f\in L_{1,loc}(D)$ and a constant $K>0,$ such that for any $(s,y,z)\in[0,T]\times D\times R^d,$ $|F_D(s,y,z)|\leq f(y)|z|^2+K|z|.$\\
\\
\textbf{Proposition 4.4}\ \ \emph{Let $F_D(s,y,z)$ satisfy Assumption (A). If $\xi\in L_D({\cal{F}}_T)$ and the range of $\xi$ is included in a closed interval $[a,b]\subset D,$ then the BSDE$(F_D(s,y,z),\xi)$ has a solution $(Y_t,Z_t)\in{\mathcal{S}}^\infty\times{\cal{H}}^2_{BMO},$ and the range of $Y$ is included in $[a,b]$.}\\\\
\emph{Proof.}  Given $\delta\in D.$ Let $\xi_1=\xi1_{\{\xi< \delta\}}+\delta1_{\{\xi\geq\delta\}}$ and $\xi_2=\xi1_{\{\xi\geq \delta\}}+\delta1_{\{\xi<\delta\}}$. By Proposition 3.3, BSDE$(-f(y)|z|^2-K|z|,\xi_1)$ and BSDE$(f(y)|z|^2+K|z|,\xi_2)$ have unique solutions $(Y^1,Z^1)\in{\mathcal{S}}^\infty\times{\cal{H}}^2_{BMO}$ and $(Y^2,Z^2)\in{\mathcal{S}}^\infty\times{\cal{H}}^2_{BMO},$ respectively, such that the range of $Y^i$ is included in $[a,b], i=1,2.$ Clearly, we have $\xi_1\leq \delta\leq\xi_2.$ Since the BSDE$(0,\delta)$ has a unique solution $(\delta,0),$ by comparison theorems (see Proposition 4.1 and Corollary 4.2), we get $Y^1\leq \delta\leq Y^2.$ Then by the method of Theorem 4.1 in \cite{B17}, we can prove that the BSDE$(F_D(s,y,z),\xi)$ has at least one solution $(Y_t,Z_t)\in{\mathcal{S}}\times{{H}}^2$ such that $Y^1\leq Y\leq Y^2,$ which implies that the range of $Y$ is included in $[a,b].$ We sketch this proof. Firstly, one can show that the reflected BSDE$(F_D(s,y,z),\xi)$ with lower obstacle $Y^1$ and upper obstacle $Y^2$ has at least one solution $(Y_t,Z_t,K_t^+,K_t^-)$ such that $(Y_t,Z_t)\in{\mathcal{S}}\times{{H}}^2$. Then we can further show $K_t^+=K_t^-=0,$ which completes this proof.

Now we prove $Z_t\in{\cal{H}}^2_{BMO}.$  We define
$$\bar{u}_f(x):=\int^x_\alpha\exp\left(2\int^y_\alpha f(z)dz\right)\left(\int^y_\alpha\exp\left(-2\int^z_\alpha f(s)ds\right)dz\right)dy,\ \ x\in D.$$
By Lemma A.1(II) in \cite{B17}, we get that $\bar{u}_f(x)$ belong to $W^2_{1,loc}(D)\cap C^{1}(D)$ and
$$\frac{1}{2}\bar{u}''_f(x)-f(x)\bar{u}'_f(x)=\frac{1}{2},\ a.e. x\in D.\eqno(32)$$
For $n\geq1$, we define the stopping time $$\tau_n=\inf\left\{t\geq0,\int_0^t|Z_s|^2ds\geq n\right\}\wedge T.$$
For $\tau\in{\cal{T}}_{0,T},$ applying Lemma 2.3 to $\bar{u}_f(Y_t)$ on $[\tau,\tau\vee\tau_n]$ and by Assumption (A), we have
\begin{eqnarray*}
\bar{u}_f(Y_\tau)&=&\bar{u}_f(Y_{\tau\vee\tau_n})+\int_\tau^{\tau\vee\tau_n}(\bar{u}_f'(Y_s)F_D(s,Y_s,Z_s)
-\frac{1}{2}\bar{u}''_f(Y_s)|Z_s|^2)ds-\int_\tau^{\tau\vee\tau_n}\bar{u}'_f(Y_s)Z_sdB_s\\
&\leq&\bar{u}_f(Y_{\tau\vee\tau_n})+\int_\tau^{\tau\vee\tau_n}(\bar{u}_f'(Y_s)( K|Z_s|+f(Y_s)|Z_s|^2)-\frac{1}{2}\bar{u}''_f(Y_s)|Z_s|^2)ds
\\&&-\int_\tau^{\tau\vee\tau_n}\bar{u}'_f(Y_s)Z_sdB_s\\
&\leq&\bar{u}_f(Y_{\tau\vee\tau_n})+\int_\tau^{\tau\vee\tau_n}(K^2|\bar{u}_f'(Y_s)|^2
+\frac{1}{4}|Z_s|^2)ds\\&&+\int_\tau^{\tau\vee\tau_n} (\bar{u}_f'(Y_s)f(Y_s)-\frac{1}{2}\bar{u}_f''(Y_s))|Z_s|^2)ds-\int_\tau^{\tau\vee\tau_n}\bar{u}_f'(Y_s)Z_sdB_s.
\end{eqnarray*}
which together with (32) gives
\begin{eqnarray*}
E\left[\int_\tau^{\tau\vee\tau_n}|Z_s|^2ds|{\cal{F}}_\tau\right]&\leq&4E\left[\bar{u}_f(Y_{\tau\vee\tau_n})+\int_\tau^{\tau\vee\tau_n} K^2|\bar{u}_f'(Y_s)|^2ds-\bar{u}_f(Y_{\tau})|{\cal{F}}_\tau\right].
\end{eqnarray*}
Since $Y$ is bounded and $\bar{u}_f(x)\in C^{1}(D),$ by Fatou's lemma, we have
\begin{eqnarray*}
E\left[\int_\tau^T|Z_s|^2ds|{\cal{F}}_\tau\right]\leq\liminf_{n\rightarrow\infty}E\left[\int_\tau^{\tau\vee\tau_n}|Z_s|^2ds|{\cal{F}}_\tau\right]<\bar{K},
\end{eqnarray*}
where $\bar{K}$ is a constant independent of $\tau.$ Thus, $Z_t\in{\cal{H}}^2_{BMO}.$ The proof is complete. \ \ $\Box$\\\\
\textbf{Proposition 4.5}\ \ \emph{Let $f_1, f_2\in L_{1,loc}(D),$ and let both be local Lipschitz continuous and satisfy a linear growth condition. For $\xi\in L_D({\cal{F}}_T)$ whose range is included in a closed subset of $D,$ let $(Y^{\xi,i}_t,Z^{\xi,i}_t)\in{\mathcal{S}}^\infty
\times{\cal{H}}^2_{BMO}$ be the solution of BSDE$(f_i(y)|z|^2,\xi), i=1,2$. If for all such $\xi$ and each $t\in[0,T],$ we have $Y^{\xi,1}_t\geq Y^{\xi,2}_t,$ then for each $a\in D,$ we have $f_1(a)\geq f_2(a).$}\\
\\
\emph{Proof.} For each $m\geq1,$ set $$\Theta_m:=\{a\in D: f_1(a)< f_2(a)-\frac{1}{m}\}.$$
Since $f_1$ and $f_2$ are both continuous, if $f_1\geq f_2$ does not hold true, then we will get that there exists a integer $n\geq1,$ such that $\lambda(\Theta_n)>0,$ where $\lambda$ is the Lebesgue measure. For $a\in \Theta_n$ and $\delta\in R^d\ (\delta\neq0)$, we consider the following two SDEs
$$Y^i_t=a-\int_0^{t}f_i(Y^i_s)|\delta|^2ds+\int_0^{t}\delta dB_s, \ \ i=1,2.\eqno(33)$$
Since $f_1$ and $f_2$ both are local Lipschitz continuous and satisfy a linear growth condition, the SDEs in (33) both have unique solutions $Y^1_t\in{\mathcal{S}}^2$ and $Y^2_t\in{\mathcal{S}}^2$, respectively. Clearly, $\Theta_n$ is an open set, which means that there exists a constant $K>0$ such that $[a-K,a+K]\subset \Theta_n.$ We define the stopping times
\begin{eqnarray*}
\tau^i&:=&\inf\left\{t\geq0: Y^i_t\notin(a-K,a+K)\right\}\wedge T, \ i=1,2;\\
\tau^3&:=&\inf\left\{t\geq0: f_1(Y^1_t)\geq f_2(Y^2_t)-\frac{1}{2n}\right\}\wedge T;\\
\tau&:=&\tau^1\wedge\tau^2\wedge\tau^3.
\end{eqnarray*}
Then we have $a-K\leq Y^i_\tau\leq a+K, i=1,2,$ $P(\tau>0)=1$ and $f_2(Y^2_s)-f_1(Y^1_s)\geq\frac{1}{2n}, s\in[0,\tau].$ Thus by (33), we have
$$Y^1_\tau-Y^2_\tau=\int_0^\tau (f_2(Y^2_s)-f_1(Y^1_s))|\delta|^2ds\geq\frac{1}{2n}|\delta|^2\tau>0.\eqno(34)$$
Set
$$y^i_t:=\left\{
      \begin{array}{ll}
        Y^i_t, & t\in[0,\tau] \\
        Y^i_\tau, & t\in(\tau,T]
      \end{array}
   \right.,\ \ \
z^i_t:=\left\{
      \begin{array}{ll}
        \delta, & t\in[0,\tau] \\
        0,& t\in(\tau,T]
      \end{array}
   \right.\ \ i=1,2.
$$
We can find that BSDE$(f_i(y)|z|^2,Y^i_\tau)$ has a solution $(y^i_t,z^i_t)\in{\mathcal{S}}^\infty
\times{\cal{H}}^2_{BMO}$ and the range of $y^i_t$ is included in $[a-K,a+K],$ $i=1,2$. Let $(y^3_t,z^3_t)\in{\mathcal{S}}^\infty
\times{\cal{H}}^2_{BMO}$ be the solution of BSDE$(f_2(y)|z|^2,Y^1_\tau)$ such that the range of $y^3_t$ is included in $[a-K,a+K]$. Then by Lemma 2.1(i), it is not hard to check that $(\int_0^tu'_{f_2}(y^i_s)z^i_sdB_s)_{t\in[0,T]}$ is a martingale and $u_{f_{2}}(y^i_t)\in{\mathcal{S}}^p$, $i=2,3$. Then by (34) and the strict comparison theorem in Proposition 4.1, we have
$$a=y^2_0<y^3_0.$$
By the assumptions of the proposition, for any $\xi\in L_D({\cal{F}}_T)$ whose range is included in a closed subset of $D,$ we have $Y^{\xi,1}_t\geq Y^{\xi,2}_t.$  Thus we have $y^3_0\leq y^1_0=a,$ which together with the above inequality induces a contradiction. The proof is complete.\ \ $\Box$
\section{Application of BSDE($f(y)|z|^2,\xi)$ to PDEs}
In this section, we will give an application to the following quadratic PDE
$$\left\{
    \begin{array}{ll}
      \partial_tv(t,x)+{\cal{L}}v(t,x)+f(v(t,x))|\sigma^*\nabla_xv(t,x)|^2=0,\ \ (t,x)\in[0,T)\times R^d,\\
      v(T,x)=g(x),\ \ x\in R^d,
    \end{array}
  \right.\eqno(35)
$$
where $f\in L_{1,loc}(D)$, $g(x):R^d\mapsto D,$ and ${\cal{L}}$ is the infinitesimal generator of the solution $X_s^{t,x}$ of SDE	
$$X_s^{t,x}=x+\int_t^sb(r,X_r^{t,x})dr+\int_t^s\sigma(r,X_r^{t,x})dB_r,\ x\in R^d,\ s\in[t,T],$$
where $b:[0,T]\times R^d\mapsto R^d,\ \sigma:[0,T]\times R^d\mapsto R^{d\times d}.$ ${\cal{L}}$ is a second order differential generator given by
$${\cal{L}}:=\frac{1}{2}\sum_{i,j=1}^d(\sigma\sigma^*)_{i,j}(s,x)\frac{\partial^2}{\partial_{x_i}\partial_{x_j}}
+\sum_{i=1}^db_{i}(s,x)\frac{\partial}{\partial_{x_i}}.$$ \\
\textbf{Assumption (B)}
(i) $f(\cdot)$ is continuous and nonnegative, and $g(\cdot)$ is measurable;

(ii) $u_f(g(\cdot))$ is continuous and has a polynomial growth;

(iii)  $b(t,\cdot)$ and $\sigma(t,\cdot)$ are uniformly Lipschitz continuous and have linear growth.\\

Let Assumption (B) hold. By Theorem 4.4 on page 61 in \cite{M}, we have for each $p\geq1$, $X^{t,x}_T\in L^p({\cal{F}}_T),$ which implies for each $p\geq1$, $u_f(g(X^{t,x}_T))\in L^p({\cal{F}}_T).$ Then by Theorem 3.2, the following Markovian BSDE
$$Y_s^{t,x}=g(X_T^{t,x})+\int_s^Tf(Y_r^{t,x})|Z_r^{t,x}|^2dr-\int_s^TZ_r^{t,x}dB_r,\ \ s\in[t,T],\eqno(36)$$
has a unique solution $(Y_s^{t,x},Z_s^{t,x})\in{\mathcal{S}}
\times{{H}}^2$ such that $u_f(Y_s^{t,x})=E[u_f(g(X^{t,x}_T))|{\cal{F}}_s].$ Set $v(t,x):=Y_t^{t,x}.$  Then we have the following Theorem 5.1.\\
\\
\textbf{Theorem 5.1}\ \ Under Assumption (B), $v(t,x)$ is a viscosity solution of PDE(35). \\\\
\emph{Proof.} By Assumption (B), (36) and Theorem 3.2(v), we get that the Markovian BSDE$(0,u_f(g(X^{t,x}_T)))$ has a unique solution $(y_s^{t,x},z_s^{t,x})=(u_f(Y_s^{t,x}),u'_f(Y_s^{t,x})Z_s^{t,x})\in{\mathcal{S}}^p
\times{\cal{H}}^p, p>1$. Since $u_f(v(t,x))=u_f(Y_t^{t,x})=y_t^{t,x},$ by Theorem 3.2 in \cite{El2}, we get that $u_f(v(t,x))$ is continuous in $(t,x)$ and grows at most polynomially at infinity. Using the uniqueness of the solution to the Markovian BSDE$(0,u_f(g(X^{t,x}_T)))$, we have $$u_f(Y^{t,x}_s)=y^{t,x}_s=y^{s,X_s^{t,x}}_s=u_f(v(s,X^{t,x}_s)), \ s\in[t,T].$$ Then by Lemma 2.1(v), we get that $v(t,x)$ is continuous in $(t,x)$ and for $s\in[t,T],$
$$Y^{t,x}_s=v(s,X^{t,x}_s).\eqno(37)$$
We will show that $v(t,x)$ is a viscosity subsolution of PDE(35). Let $\phi(t,x)\in C^{1,2}([0,T]\times R^d)$ and $(t,x)$ be a local maximum point of $v-\phi$. Without loss of generality, we assume $v(t,x)=\phi(t,x).$ If
$$\frac{\partial\phi}{\partial t}(t,x)+{\cal{L}}\phi(t,x)+f(v(t,x))|\sigma^*\nabla_x\phi(t,x)|^2<0,$$
then by continuity, there exist $\beta\in(0,T-t],$ $c>0$ and $K>0,$ such that for $s\in[t,t+\beta]$ and $y\in[x-c,x+c]$, we have
$$v(s,y)\leq\phi(s,y),\ \ \textmd{and} \ \ \frac{\partial\phi}{\partial t}(s,y)+{\cal{L}}\phi(s,y)+f(v(s,y))|\sigma^*\nabla_x\phi(s,y)|^2\leq-K.\eqno(38)$$
We define a stopping time $\tau=\inf\{s\geq t; |X_s^{t,x}-x|\geq c\}\wedge(t+\beta),$ then $t<\tau \leq t+\beta$ and $x-c\leq X_{s\wedge\tau}^{t,x}\leq x+c,  s\in[t,t+\beta].$ By (36) and (37), we deduce that the following BSDE
$$\bar{Y}_s=v(\tau,X_\tau^{t,x})+\int_s^{t+\beta} 1_{\{r\leq\tau\}}f(v(r,X_r^{t,x}))|\bar{Z}_r|^2dr-\int_s^{t+\beta}\bar{Z}_rB_r,\ \ s\in[t,t+\beta],\eqno(39)$$
has a solution $(\bar{Y}_s,\bar{Z}_s)=({Y}_{s\wedge\tau}^{t,x},1_{\{s\leq\tau\}}{Z}_s^{t,x})\in{\mathcal{S}}^\infty
\times{{H}}^2.$ Applying It\^{o}'s formula to $\phi(s,X_s^{t,x})$ on $s\in[t,(t+\beta)\wedge\tau]$, we deduce that the following BSDE
$${\tilde{Y}}_s=\phi(\tau,X_\tau^{t,x})+\int_t^{t+\beta} -1_{\{r\leq\tau\}}\left(\frac{\partial\phi}{\partial t}(r,X_r^{t,x})+{\cal{L}}\phi(r,X_r^{t,x})\right)dr-\int_t^{t+\beta}{\tilde{Z}}_rB_r,\ \ s\in[t,t+\beta],\eqno(40)$$
has a solution $({\tilde{Y}}_s,{\tilde{Z}}_s)=(\phi({s\wedge\tau},X_{s\wedge\tau}^{t,x}), 1_{\{s\leq\tau\}}\sigma^*\nabla_x\phi(s,X_s^{t,x}))\in{\mathcal{S}}^\infty
\times H^2.$ By (38), we have $v(\tau,X_\tau^{t,x})\leq\phi(\tau,X_\tau^{t,x})$ and
$$\int_t^{t+\beta}1_{\{r\leq\tau\}}\left(-\frac{\partial\phi}{\partial t}(s,X_r^{t,x})-{\cal{L}}\phi(s,X_r^{t,x})-f(v(r,X_r^{t,x}))|\sigma^*\nabla_x\phi(r,X_r^{t,x})|^2)\right)dr\geq K(\tau-t)>0.$$
By Assumption (B)(i), we get that $\gamma_r:=1_{\{r\leq\tau\}}f(v(r,X_r^{t,x}))$ is nonnegative and bounded in $[t,t+\beta].$ Thus (39) is indeed a BSDE with generator $\gamma_r|z|^2,$ which is convex in $z$. Then by (39), (40) and the strict comparison theorem (see Theorem 5 in \cite{BH08}), we have $v(t,x)=\bar{Y}_t<\tilde{Y}_t=\phi(t,x).$ This induces a contradiction. Thus $v(t,x)$ is a viscosity subsolution of PDE(35). Similarly, we can also show that $v(t,x)$ is a viscosity supersolution of PDE(35). The proof is complete.\ \ $\Box$\\\\
\textbf{Remark 5.1}  \begin{itemize}
                     \item The viscosity solution of PDE(35) with $f(y)=\frac{1}{y}$ was considered using an approximating method in Theorem 4.5 in \cite{BT}. The corresponding weak solution problems can be found in \cite{D}, where the initial condition corresponding to the terminal condition $g(x)$ in PDE(35), is assumed to be bounded. It is not hard to find that if the range of $g(x)$ is included in a closed subset of $D$, then Theorem 5.1 will hold true for any continuous nonnegative function $f\in L_{1,loc}(D).$
                     \item In Theorem 5.1, $f(x)$ is assumed to be nonnegative in order to use the convexity in the strict comparison theorem (see Theorem 5 in \cite{BH08}). Correspondingly, when $f(x)$ is nonpositive, the concavity will be obtained. Then by the strict comparison theorem, Theorem 5.1 will also hold true.
                     \item In Theorem 5.1, we only obtain the existence of a viscosity solution. The uniqueness of viscosity solutions to PDEs is usually obtained by using a comparison theorem (see \cite{DL} and \cite{K}, etc), but we still cannot establish a comparison result in our case. A uniqueness result for quadratic PDEs can be found in \cite{DL} or the Remark on page 566 in \cite{BH08}. In some cases, the viscosity solution in Theorem 5.1 has a polynomial growth. We show an example. Let $D=[a,\infty),\ a\in R.$ Since $f\geq0$, by Lemma 2.1(vi), we have $u_f(x)\geq x-\alpha,$ which together with Lemma 2.1(ii) gives $u_f^{-1}(x)\leq x+\alpha.$ Since $u_f(v(t,x))=y_t^{t,x},$ by Theorem 3.2 in \cite{El2}, we have $|u_f(v(t,x))|\leq K(1+|x|^p),$ for some $p>1.$ Then by the fact that $u_f^{-1}$ is strictly increasing, we deduce $$a\leq v(t,x)\leq K(1+|x|^p)+\alpha,$$
                         which implies that $v(t,x)$ has a polynomial growth.
                   \end{itemize}\
\\
\textbf{Acknowledgements}\ \   The authors would like to thank the anonymous referee for their valuable comments and suggestions. The first author is supported by the National Natural Science Foundation of China (No. 11571024) and the Humanities and Social Science Foundation of Ministry of Education of China (No. 20YJCZH245).


\begin{thebibliography}{99}
\addtolength{\itemsep}{-0.8 em}
\bibitem{BE} Briand P. Elie R. A simple constructive approach to quadratic BSDEs with or without delay. Stochastic Processes and their Applications, 123, 2921-2939 (2012)
\bibitem{B19} Bahlali K. Solving unbounded quadratic BSDEs by a domination method. arXiv:1903.11325.(2019)
\bibitem{B17} Bahlali K., Eddahbi M., and Ouknine Y. Quadratic BSDEs with $L^2$-terminal data
existence results, Krylov's estimate and It\^{o}-Krylov's formula. Ann. Probab. 45(4), 2377-2397 (2017)
\bibitem{BT} Bahlali K., Tangpi L. BSDEs driven by $|z|^ 2/y $ and applications. arXiv:1810.05664v2 (2018)
\bibitem{lp}  Briand P., Delyon B. Hu Y., Pardoux E. and Stoica L. $L^p$ of backward stochastic differential equations. Stochastic processes and their applications, 108, 109-129 (2003)
\bibitem{BH06}  Briand P., and Hu Y. BSDE with quadratic growth and unbounded terminal
value. Probability Theory and Related Fields, 136(4), 604-618 (2006)
\bibitem{BH08}  Briand P., and Hu Y. Quadratic BSDEs with convex generators and unbounded terminal
value. Probability Theory and Related Fields, 141, 543-567 (2008)
\bibitem{D}  Dall'Aglio A., Orsina L., and Petitta F. Existence of solutions of degenerate parabolic equations with singular terms. Nonlinear
Analysis, 131, 273-288 (2016)
\bibitem{DL} Da Lio, F., Ley, O., Uniqueness results for second-order Bellman-Isaacs equations under quadratic
growth assumptions and applications. SIAM J. Control Optim. 45(1), 74-106 (2006)
\bibitem{El} El Karoui N. and Rouge R. Pricing via utility maximization and entropy. Math. Finance 10(2), 259-276 (2000)
\bibitem{El2} El Karoui N., Hamad\`{e}ne S., and Matoussi A. Backward stochastic differential equations and applications. In R. Carmona, editor, Indifference Pricing: Theory and Applications, Princeton Series in Financial Engineering, pages 267-320. Princeton University Press, (2009)
\bibitem{H} Hu Y., Imkeller P., M\"{u}ller M. Utility maximization in incomplete markets. Ann. Appl Probab. 15(3), 1691-1712 (2005)
\bibitem{ks} Karatzas I., Shreve S.
Brownian Motion and Stochastic Calculus,
Second Edition. Springer (1991)
\bibitem{K}  Kobylanski M. Backward stochastic differential equations and partial differential equations with
quadratic growth. Ann. Probab. 28(2), 558-602 (2000)
\bibitem{M}  Mao X. Stochastic differential
equation with applications, Second edition.  Woodhead Publishing  (2007)
\bibitem{PP}  Pardoux E. and Peng, S. Adapted solution of a backward stochastic differential
equation. Systems Control Lett. 14 55-61 (1990)
\bibitem{Pr}  Protter P. Stochastic integration and differential
equation. Springer (2005).
\bibitem{Y}  Yang H. $L^p$-solutions of quadratic BSDEs. arXiv:1506.08146v2 (2017)
\end{thebibliography}
\end{document}